\documentclass[12pt]{amsart} 
\usepackage{amsmath} \usepackage{amssymb} \usepackage{amsthm}
\usepackage{amscd} \usepackage[all]{xy} 
\def\a{{\mathfrak{a}}} \def\b{{\mathfrak{b}}}
 \def\J{{\mathcal{J}}} \def\m{{\mathfrak{m}}} \def\Z{{\mathbb{Z}}}
 \def\O{{\mathcal{O}}}
\def\Q{{\mathbb{Q}}} \def\R{{\mathbb{R}}} \def\C{{\mathbb{C}}}
 \def\Ker{{\mathrm{Ker\;}}} \def\hight{{\mathrm{ht\; }}}
  
\def\Spec{{\mathrm{Spec\; }}} 
\def\Div{{\mathrm{div}}}
\theoremstyle{plain}
\newtheorem{thm}{Theorem}[section] 
\newtheorem{prop}[thm]{Proposition}
\newtheorem{conj}[thm]{Conjecture}
 
\newtheorem{lem}[thm]{Lemma}
\theoremstyle{definition} 
\newtheorem{dfn}[thm]{Definition}
\newtheorem{eg}[thm]{Example} 
\theoremstyle{remark}
\newtheorem{rmk}[thm]{Remark}
\newtheorem*{cl}{Claim}

\title{On F-pure thresholds}
\author{Shunsuke Takagi}
\address{Faculty of Mathematics, Kyushu University, 
6-10-1 Hakozaki, Higashi, Fukuoka, 812-8581 JAPAN}
\email{stakagi@math.kyushu-u.ac.jp}
\author{Kei-ichi Watanabe}
\address{Department of Mathematics, College of Humanities and Sciences, 
Nihon University, Setagaya-ku, Tokyo 156--0045, JAPAN}
\email{watanabe@math.chs.nihon-u.ac.jp}
\subjclass[2000]{13A35, 14B05}
\dedicatory{Dedicated to Professor~Melvin~Hochster on the~occasion of his~sixtieth~birthday.}
\begin{document}

\maketitle
\markboth{SHUNSUKE TAKAGI AND KEI-ICHI WATANABE}{ON F-PURE THRESHOLDS}

\begin{abstract}
Using the Frobenius map, we introduce a new invariant for a pair $(R,\a)$ of a ring $R$ and an ideal $\a \subset R$, which we call the F-pure threshold $\mathrm{c}(\a)$ of $\a$, and study its properties. We see that the F-pure threshold characterizes several ring theoretic properties. 
By virtue of Hara and Yoshida's result \cite{HY}, the F-pure threshold $\mathrm{c}(\a)$ in characteristic zero corresponds to the log canonical threshold $\mathrm{lc}(\a)$ which is an important invariant in birational geometry. 
Using the F-pure threshold, we prove some ring theoretic properties of three-dimensional terminal singularities of characteristic zero. 
Also, in fixed prime characteristic, we establish several properties of the F-pure threshold similar to those of the log canonical threshold with quite simple proofs. 
\end{abstract}

\section*{Introduction}
The log canonical threshold is an important invariant in higher dimensional birational geometry; it measures how far a pair $(X, Y)$ is from being log canonical where $X$ is an algebraic variety and $Y$ is a closed subscheme of $X$. 
This invariant is defined in terms of resolution of singularities and discrepancy divisors, and it is in general quite mysterious. 
In order to understand its behavior, it is helpful to formulate the log canonical threshold purely algebraically. 

In \cite{T}, the first-named author generalized the notions of F-pure and F-regular rings to those for a pair $(R, \a^t)$ of a ring $R$ of characteristic $p>0$ and an ideal $\a \subset R$ with real exponent $t\ge 0$. 
When $R$ is an F-pure ring, the pair $(R,\a^t)$ will be F-pure for sufficiently small $t \ge 0$, whereas it cannot be so when $t$ is very large. 
We define the F-pure threshold $\mathrm{c}(\a)$ of $\a$ to be the critical value of $t$ for this property. % (Definition \ref{def}). 
We can also define the F-pure threshold for an ideal of a ring of characteristic zero, using reduction to characteristic $p>0$. % (Definition \ref{char 0}). 
Then, by virtue of Hara and Yoshida's theorem \cite[Theorem 6.8]{HY}, 
%the F-pure threshold corresponds to the log canonical threshold. Precisely speaking, for a $\Q$-Gorenstein normal local ring $R$ of characteristic zero, in the case where $R$ is of strongly F-regular type (or equivalently a log terminal singularity), 
the F-pure threshold in characteristic zero coincides with the log canonical threshold under mild conditions. 
%when $R$ is a $\Q$-Gorenstein normal local ring of strongly F-regular type (or equivalently a log terminal singularity). 
We expect that this coincidence holds true whenever both invariants exist, and that the F-pure threshold might shed new light on the theory of the log canonical threshold. 

In this paper, we study properties of the F-pure threshold. 
Although our first motivation was to investigate the log canonical threshold via the F-pure threshold, we find that the F-pure threshold itself is an interesting invariant in commutative algebra. 
For example, when $R$ is a Gorenstein strongly F-regular graded ring generated by elements of degree one over a field, the F-pure threshold $\mathrm{c}(R_+)$ of the irrelevant ideal of $R$ is equal to the minus of a-invariant, which is a classical invariant in commutative algebra introduced by Goto and the second-named author \cite{GW}. 

In section $2$, we study fundamental properties of the F-pure threshold in characteristic $p>0$. 
Specifically, for a $d$-dimensional Noetherian local ring $(R,\m)$ of characteristic $p>0$, we prove that if the F-pure threshold $\mathrm{c}(\m)$ of the maximal ideal $\m$ of $R$ is greater than $d-r$, then $\m^r$ is contained in any minimal reduction $J \subset \m$ (Proposition \ref{key lemma} (2)). In the case where $R$ is $\Q$-Gorenstein Cohen-Macaulay, we have a more powerful result (Proposition \ref{reduction}) of this. 
Also, we see that $\mathrm{c}(\m)$ characterizes several ring theoretic properties: $\mathrm{c}(\m)>d-1$ if and only if $R$ is regular, and when $R$ is $\Q$-Gorenstein Cohen-Macaulay, $\mathrm{c}(\m)=d-1$ if and only if the $\m$-adic completion $\widehat{R}$ of $R$ is the generalized $\mathrm{cA}_n$-singularity. (Theorem \ref{big value}). 

In section $3$, we investigate properties of the F-pure threshold in characteristic zero, especially the relationship with three-dimensional terminal singularities. 
%Although terminal singularities are geometric objects, 
%we can prove some ring theoretic properties for three-dimensional terminal singularities using the F-pure threshold.
we give a characterization of three-dimensional Gorenstein terminal singularities in terms of the F-pure threshold (Proposition \ref{Gorenstein terminal}), and generalize Kakimi's result \cite{Ka} on the multiplicity of non-Gorenstein three-dimensional terminal singularities (Proposition \ref{kakimi}). 

In section $4$, we reprove or translate some results related to the log canonical threshold using the F-pure threshold in characteristic $p>0$. 
For a regular local ring $(R,\m)$ of characteristic $p>0$, we prove that the F-pure threshold $\mathrm{c}(\a)$ of an $\m$-primary ideal $\a \subset R$ gives a lower bound for the multiplicity of $\a$ by an argument similar to that in \cite{DEM} (Proposition \ref{multiplicity}). 
We also show the restriction property and the summation property of the F-pure threshold by simple ring theoretic arguments (Proposition %\ref{bound}, \ref{bound2}, 
\ref{restriction} and \ref{subadditivity}), whereas corresponding results on the log canonical threshold need a deep vanishing theorem (see \cite{DK} or \cite{Ko}) or a characterization of the log canonical threshold via jet schemes due to Musta{\c{t}}{\v{a}} \cite{Mu}. 

{\bf Acknowledgments.} The first-named author is grateful to Melvin Hochster and Robert Lazarsfeld for their hospitality during his stay at the University of Michigan in the fall of 2003. 
The authors are indebted to Mircea Musta{\c{t}}{\v{a}} for pointing out some mistakes in a previous version of this paper. The first-named author was partially supported by the Japan Society for the Promotion of Science Research Fellowships for Young Scientists. 
The second-named author was supported in part by Grant aid in Scientific Researches, \# 13440015 and \# 13874006.

\section{Preliminaries}
\subsection{Log canonical threshold}
First we recall the definitions of singularities of pairs and the log canonical threshold. Refer to \cite{Ko} and \cite{KM} for the detail.

Let $X$ be a $\Q$-Gorenstein normal variety over a field of characteristic zero, $Y \subseteq X$ a closed subscheme defined by an ideal sheaf $\a \subset \O_X$ and $t \ge 0$ a real number.
Suppose $f:\widetilde{X} \to X$ is a log resolution of the pair $(X, Y)$. 
 That is, $f$ is a proper birational morphism with $\widetilde{X}$ 
 nonsingular such that the ideal sheaf $\a\O_{\widetilde{X}}=\O_{\widetilde{X}}(-F)$ is invertible 
and  $\mathrm{Supp \;} F \cup \mathrm{Exc}(f)$ is a simple normal 
crossing divisor, where $\mathrm{Exc}(f)$ is the exceptional locus of $f$ 
(Hironaka \cite{Hi} proved that log resolutions always exist). 
Let $K_X$ and $K_{\widetilde{X}}$ denote canonical divisors of $X$ and 
$\widetilde{X}$, respectively. 
Then there are finitely many irreducible (not necessarily exceptional) divisors $E_i$ on $\widetilde{X}$ and real numbers $a_i$ so that there exists an $\R$-linear equivalence of $\R$-divisors 
%$\sum_{i} a_iE_i+tF$ is an exceptional divisor satisfying 
$$K_{\widetilde{X}} \underset{\text{$\Q$-lin.}}{\sim} f^*K_X+
\sum_{i} a_iE_i+tF.$$

\begin{dfn}\label{pair}
Under the notation as above:
\renewcommand{\labelenumi}{(\roman{enumi})}
\begin{enumerate}
\item 
We say that the pair $(X,tY)$ is \textit{terminal} if $a_i >0$ for all $i$. 
\item We say that the pair $(X,tY)$ is \textit{Kawamata log terminal} (or \textit{klt} for short) if $a_i>-1$ for all $i$. 
\item We say that the pair $(X,tY)$ is \textit{log canonical} (or \textit{lc} for short) if $a_i \ge -1$ for all $i$. 
\item 
Suppose that $X$ has at worst log canonical singularities at a point $x \in X$. Then 
we define the \textit{log canonical threshold} (or \textit{lc threshold} for short) of $\a$ at $x$ to be
$$\mathrm{lc}(\a;x)=\sup \{s \in \R_{\ge 0} \mid \textup{the pair $(X,sY)$ is lc at $x$} \}. $$
\item 
The \textit{multiplier ideal sheaf} $\mathcal{J}(X, \a^t)$ associated to $\a$ with exponent $t$ is defined to be the ideal sheaf
$$\mathcal{J}(X, \a^t)=f_*\mathcal{O}_{\widetilde{X}}(\sum_i \lceil a_i \rceil E_i) \subset \O_X.$$
\end{enumerate}
These definitions are independent of the choice of a log resolution $f:\widetilde{X} \to X$.
\end{dfn}

\begin{rmk}\label{multiplier ideals}
$(1)$ 
$X$ has at worst terminal (resp. log terminal, log canonical) singularities if and only if the pair $(X,0)$ is terminal (resp. klt, lc).

$(2)$
The pair $(X, tY)$ is klt if and only if $\mathcal{J}(X, \a^t)=\O_X$.

$(3)$
Let $\mathcal{I} \subset \O_X$ be the radical ideal sheaf which defines the non-log-terminal locus of $X$. 
If $\mathcal{I} \supset \a$, then
$$\mathrm{lc}(\a)=\sup \{s \in \R_{\ge 0} \mid \textup{the pair $(X, sY)$ is klt}\}.$$

$(4)$
Set $c=\mathrm{lc}(\a)$, then the pair $(X,cY)$ is lc. 
\end{rmk}

\subsection{F-singularities of pairs}
Here we briefly review the definition and fundamental properties of F-singularities of pairs which we need later. 
Refer to \cite{HY} and \cite{T} for the proofs. 

In this paper, all rings are reduced commutative rings with unity. 
For a ring $R$, we denote by $R^{\circ}$ the set of elements of $R$ which are not in any minimal prime ideal. 
Let $R$ be a ring of characteristic $p >0$ and $F\colon R \to R$ the Frobenius map which sends $x \in R$ to $x^p \in R$. 
For an integer $e \ge 0$, the ring $R$ viewed as an $R$-module via the $e$-times iterated Frobenius map $F^e \colon R \to R$ is denoted by ${}^e\! R$. 
Since $R$ is assumed to be reduced, we can identify $F^e \colon R \to {}^e\! R$ with the natural inclusion map $R \hookrightarrow R^{1/p^e}$. We say that $R$ is {\it F-finite} if ${}^1
\! R$ (or $R^{1/p}$) is a finitely generated $R$-module. 
We assume always that all rings of characteristic $p>0$ are F-finite throughout this paper. 

\begin{dfn}\label{F-def}
Let $\a$ be an ideal of an F-finite reduced ring $R$ of characteristic $p>0$ such that $\a \cap R^{\circ} \ne \emptyset$, and let $t \ge 0$ be a real number. 
\renewcommand{\labelenumi}{(\roman{enumi})}
\begin{enumerate}
\item
The pair $(R, \a^t)$ is said to be {\it F-pure} (or $R$ is said to be {\it F-pure} with respect to $\a$ and $t$) if for all large $q=p^e \gg 0$, there exists an element $d \in \a^{\lfloor t(q-1) \rfloor}$ such that $d^{1/q}R \hookrightarrow R^{1/q}$ splits as an $R$-module homomorphism. 
\item
The pair $(R, \a^t)$ is said to be {\it strongly F-pure} (or $R$ is said to be {\it strongly F-pure} with respect to $\a$ and $t$) if there exist $q=p^e>0$ and $d \in \a^{\lceil tq \rceil}$ such that $d^{1/q}R \hookrightarrow R^{1/q}$ splits as an $R$-module homomorphism.
\item
The pair $(R, \a^t)$ is said to be {\it strongly F-regular} (or $R$ is said to be {\it strongly F-regular} with respect to $\a$ and $t$) if for every $c \in R^{\circ}$, there exist $q=p^e>0$ and $d \in \a^{\lceil tq \rceil}$ such that $(cd)^{1/q}R \hookrightarrow R^{1/q}$ splits as an $R$-module homomorphism.
\end{enumerate}
\end{dfn}

\begin{rmk}
A ring $R$ is F-pure (resp. strongly F-regular) if and only if $R$ is F-pure (resp. strongly F-regular) with respect to the unit ideal $R$. 
Refer to \cite{HR}, \cite{HH1} and \cite{HH2} for properties of F-pure and strongly F-regular rings.
\end{rmk}

\begin{prop}\label{basic1}
Let notations be as in Definition \ref{F-def}.
\begin{enumerate}
\item $($\textup{cf. \cite[Proposition 3.3 (2)]{T}}$)$ 
A strongly F-regular pair is strongly F-pure, and a strongly F-pure pair is F-pure. 
\item $($\textup{\cite[Lemma 3.6]{T}}$)$
Let $c \in R^{\circ}$ be an element such that the localization $R_c$ with respect to $c$ is strongly F-regular. Then, for any ideal $\a \subseteq R$ with $\a \cap R^{\circ} \ne \emptyset$ and for any real number $t \ge 0$, the pair $(R, \a^t)$ is strongly F-regular if and only if there exist $q=p^e>0$ and $d \in \a^{\lceil tq \rceil}$ such that $(cd)^{1/q}R \hookrightarrow R^{1/q}$ splits as an $R$-linear map.
\end{enumerate}
\end{prop}

By the Matlis duality, F-singularities of pairs can be characterized via the injectivity of the induced Frobenius map on the injective hull of the residue field. 
\begin{lem}[\textup{\cite[Lemma 3.4]{T}}]\label{injective}
Let $(R,\m)$ be an F-finite reduced local ring of characteristic $p>0$, $\a \subset R$ an ideal and $t \ge 0$ a real number. 
We denote by $E_R$ the injective hull of the residue field $R/\m$ and by $F^e_R:E_R \to E_R \otimes_R {}^e R$ the induced $e$-times iterated Frobenius map on $E_R$. 
\begin{enumerate}
\item
The pair $(R,\a^t)$ is F-pure if and only if for all large $q=p^e \gg 0$, there exists an element $d \in \a^{\lfloor t(q-1) \rfloor}$ such that $dF^e:E_R \to E_R \otimes_R {}^e R$ is injective.
\item 
The pair $(R, \a^{t})$ is strongly F-pure if and only if there exist $q=p^e>0$ and $d \in \a^{\lceil tq \rceil}$ such that $dF^e:E_R \to E_R \otimes_R {}^e R$ is injective. 
\item
The pair $(R, \a^{t})$ is strongly F-regular if and only if for every element $c \in R^{\circ}$, there exist $q=p^e>0$ and $d \in \a^{\lceil tq \rceil}$ such that $cdF^e:E_R \to E_R \otimes_R {}^e R$ is injective. 

\end{enumerate}
\end{lem}

\begin{rmk}\label{split}
Let $(R,\m)$ be an F-finite reduced local ring of characteristic $p>0$ and $a, b$ elements of $R$. If ($a+b)^{1/q}R \hookrightarrow R^{1/q}$ splits as an $R$-linear map, then either $a^{1/q}R \hookrightarrow R^{1/q}$ or $b^{1/q}R \hookrightarrow R^{1/q}$ splits as an $R$-linear map. 
(Assume to the contrary that neither map splits, that is, neither of $aF^e: E_R \to E_R \otimes_R {}^e R$ nor $bF^e: E_R \to E_R \otimes_R {}^e R$ is injective.  
Then, by the definition of $E_R$,  $\Ker (a+b)F^e \supset \Ker aF^e \cap \Ker bF^e \neq 0$, which implies that $(a+b)^{1/q}R \hookrightarrow R^{1/q}$ does not split. This is a contradiction.)
In particular, let $\a \subset R$ be an ideal, $t \ge 0$ a real number and $x_1, \dots, x_m$ a system of generators of $\a$. 
Then the pair $(R,\a^t)$ is F-pure if and only if for all large $q=p^e \gg 0$, there exist nonnegative integers $l_1, \dots, l_m$ with $l_1+ \dots l_m=\lfloor t(q-1) \rfloor$ such that $(x_1^{l_1} \cdots x_m^{l_m})^{1/q}R \hookrightarrow R^{1/q}$ splits as an $R$-module homomorphism. For the other cases, we have similar characterizations. 
\end{rmk}

When the ring is a quotient of a regular local ring, we have a useful criterion for F-purity and strong F-regularity, which is called the Fedder type criterion. 
\begin{lem}[\textup{\cite[Lemma 3.9]{T}}]\label{Fedder}
Let $(R,\m)$ be an F-finite regular local ring of characteristic $p>0$ and $I \subset R$ an ideal. Fix any ideal $\a \subset R$ and any real number $t \ge 0$. Write $S=R/I$.
\begin{enumerate}
\item
The pair $(S, (\a S)^t)$ is F-pure if and only if for all large $q=p^e \gg 0$, $\a^{\lfloor t(q-1) \rfloor}(I^{[q]}:I) \not\subset \m^{[q]}$.
\item
The pair $(S, (\a S)^t)$ is strongly F-pure if and only if there exists $q=p^e>0$ such that $\a^{\lceil tq \rceil}(I^{[q]}:I) \not\subset \m^{[q]}$.
\item
The pair $(S, (\a S)^t)$ is strongly F-regular if and only if for every element $c \in R \setminus I$, there exists $q=p^e>0$ such that $c\a^{\lceil tq \rceil}(I^{[q]}:I) \not\subset \m^{[q]}$.
By Proposition $\ref{basic1}$ $(2)$, this is equivalent to saying that for an element $c \in R \setminus I$ such that the localization $S_c$ with respect to $c$ is strongly F-regular, there exists $q=p^e>0$ such that $c\a^{\lceil tq \rceil}(I^{[q]}:I) \not\subset \m^{[q]}$.
\end{enumerate}
\end{lem}

The notions of F-purity and F-regularity are also defined for a ring of characteristic zero. 
\begin{dfn}[\textup{\cite[Definition 3.7]{T}}]
Let $R$ be a reduced algebra essentially of finite type over a field $k$ of characteristic zero, $\a \subset R$ an ideal, and $t\ge0$ a real number. 
The pair $(R, \a^{t})$ is said to be of {\it dense F-pure type} (resp. {\it strongly F-regular type}) if there exist a finitely generated $\Z$-subalgebra $A$ of $k$ and a reduced subalgebra $R_A$ of $R$ essentially of finite type over $A$ which satisfy the following conditions:
\renewcommand{\labelenumi}{(\roman{enumi})}
\begin{enumerate} 
\item $R_A$ is flat over $A$, $R_A \otimes_A k \cong R$ and ${\a_A}R= \a$ where ${\a_A}=\a \cap R_A \subset R_A$. 
\item The pair $(R_{\kappa}, {\a_{\kappa}}^{t})$ is F-pure (resp. strongly F-regular) for every closed point $s$ in a dense subset of $\Spec A$, where $\kappa=\kappa(s)$ denotes the residue field of $s \in \Spec A$, $R_{\kappa}=R_A \otimes_A \kappa(s)$ and ${\a_{\kappa}}={\a_A} R_{\kappa} \subset R_{\kappa}$. 
\end{enumerate}
\end{dfn}

F-singularities of pairs in characteristic zero correspond to singularities of pairs which arise in birational geometry. 
\begin{prop}\label{klt}
Let $(R, \m)$ be a $\Q$-Gorenstein normal local ring essentially of finite type over a field of characteristic zero and write $X=\Spec R$. 
Let $Y \subsetneq X$ be a closed subscheme defined by a nonzero ideal $\a \subset R$ and $t \ge 0$ a real number.
\begin{enumerate}
\item $(\textup{\cite[Proposition 3.8]{T}})$ If the pair $(R, \a^{t})$ is of dense F-pure type, then the pair $(X, tY)$ is lc. 
\item $(\textup{\cite[Theorem 6.8]{HY}, \cite[Corollary 3.5]{T}})$ The pair $(R, \a^{t})$ is of strongly F-regular type if and only if the pair $(X, tY)$ is klt. 
\end{enumerate}
\end{prop}

\section{F-pure threshold in characteristic $p>0$}

Using the notion of F-purity for a pair $(R, \a)$ of a ring $R$ of characteristic $p>0$ and an ideal $\a \subset R$, we define the F-pure threshold. 

\begin{dfn}\label{def}
Let $R$ be a reduced F-finite F-pure ring of characteristic $p>0$ and $\a \subset R$ an ideal such that $\a \cap R^{\circ} \neq \emptyset$. 
Then we define the {\it F-pure threshold} $\mathrm{c}(\a)$ of $\a$ to be
$$\mathrm{c}(\a)=\mathrm{c}(R,\a)=\sup \{s \in \R_{\ge 0} \mid \textup{the pair $(R,\a^s)$ is F-pure} \}.$$
\end{dfn}

We collect some fundamental properties of the F-pure threshold.
\begin{prop}\label{basic}
Let notations be the same as in Definition \ref{def}. Let $\b$ be an ideal of $R$. 
\begin{enumerate}
\item If $\a \subset \b$, then $\mathrm{c}(\a) \le \mathrm{c}(\mathfrak{b})$.
\item For any integer $n>0$, $\mathrm{c}(\a)=n \cdot \mathrm{c}(\a^n)$. 
\item $\mathrm{c}(\a)=\sup \{s \in R_{\ge 0} \mid \textup{the pair $(R, \a^s)$ is strongly F-pure}\}$.
\item If the pair $(R,\a^t)$ is strongly F-regular for a real number $t\ge 0$, then $\mathrm{c}(\a)>t$. 
\item Let $I \subset R$ be the radical ideal defining the non-strongly F-regular locus of $\Spec R$. 
If $I \supset \a$, then
$$\mathrm{c}(\a)=\sup \{s \in \R_{\ge 0} \mid \textup{the pair $(R,\a^s)$ is strongly F-regular} \}.$$
\item $\mathrm{c}(\a)=\mathrm{c}(\overline{\a})$, where $\overline{\a}$ denotes the integral closure of $\a$. 
\end{enumerate}
\end{prop}
\begin{proof}
$(1)$ Obvious.

$(2)$ Since $\a^{\lfloor t(q-1) \rfloor} \subseteq (\a^n)^{\lfloor \frac{t}{n}(q-1) \rfloor}$ for all $q=p^e>0$ and $t>0$, we have $\mathrm{c}(\a) \le n \cdot \mathrm{c}(\a^n)$.
To show the converse implication, it suffices to see that if the pair $(R, (\a^n)^t)$ is F-pure for some $t>0$, then the pair $(R,\a^{nt-\epsilon})$ is also F-pure for every $nt>\epsilon>0$. 
It is, however, immediately checked by taking sufficiently large $q=p^e \gg 0$ so that $n \lfloor t(q-1) \rfloor \ge \lfloor (nt-\epsilon)(q-1) \rfloor$.

$(3)$
It is sufficient to see that if the pair $(R,\a^t)$ is F-pure for some $t>0$, then the pair $(R,\a^{t-\epsilon})$ is strongly F-pure for all $t>\epsilon>0$. It is, however, easily checked by taking sufficiently large $q=p^e \gg 0$ so that $\lfloor t(q-1) \rfloor \ge \lceil (t-\epsilon)q \rceil$.

$(4)$ 
It is enough to show that if the pair $(R,\a^t)$ is strongly F-regular for some $t>0$, then there exists a real number $1 \gg \epsilon>0$ such that the pair $(R,\a^{t+\epsilon})$ is also strongly F-regular. 
Assume that the pair $(R,\a^t)$ is strongly F-regular, and fix an arbitrary element $a \in \a \cap R^{\circ}$. 
By definition, there exist $q'=p^{e'}>0$ and $d \in \a^{\lceil tq' \rceil}$ such that $(ad)^{1/q'}R \hookrightarrow R^{1/q'}$ splits as an $R$-module homomorphism. 
Since $R$ is F-pure, the map $R^{1/q'} \hookrightarrow R^{1/qq'}$ splits for all $q=p^e>0$. 
Hence $R \xrightarrow{(ad)^{q/qq'}} R^{1/qq'}$ splits for all $q=p^e>0$.
By setting $\epsilon=1/q'>0$, we have $(ad)^q \in \a^{(1+\lceil tq' \rceil)q} \subset \a^{\lceil (t+\epsilon)qq' \rceil}$. 
Since $R$ is strongly F-regular, by Proposition \ref{basic1} $(2)$, the pair $(R,\a^{t+\epsilon})$ is strongly F-regular. 

$(5)$ 
It is obvious that 
$\mathrm{c}(\a) \ge \sup \{s \in \R_{\ge 0} \mid \textup{the pair $(R,\a^s)$ is strongly F-regular} \}$. 
To show the converse implication, it suffices to show that if the pair $(R, \a^t)$ is F-pure for some $t>0$, then the pair $(R,\a^{t-\epsilon})$ is strongly F-regular for every $t>\epsilon>0$. 
Since $I \supset \a$, by Proposition \ref{basic1} $(2)$, the pair $(R, \a^{t-\epsilon})$ is strongly F-regular if and only if for some $c \in \a \cap R^{\circ}$, there exist $q=p^e>0$ and $d \in \a^{\lceil (t-\epsilon)q \rceil}$ such that $(cd)^{1/q}R \hookrightarrow R^{1/q}$ splits as an $R$-module homomorphism. 
Fix sufficiently large $q'=p^{e'} \gg 0$ so that $\lceil (t-\epsilon)q' \rceil+1 \le \lfloor t(q'-1) \rfloor$, then the F-purity of the pair $(R, \a^t)$ implies that there exist $c \in \a \cap R^{\circ}$ and $d' \in \a^{\lceil (t-\epsilon)q' \rceil}$ such that $(cd')^{1/q'}R \hookrightarrow R^{1/q'}$ splits as an $R$-linear map (see Remark \ref{split}). 
Hence, the pair $(R, \a^{t-\epsilon})$ is strongly F-regular. 

$(6)$ 
By $(1)$, the inequality $\mathrm{c}(\a) \le \mathrm{c}(\overline{\a})$ is clear. 
So we will prove the converse inequality $\mathrm{c}(\a) \ge \mathrm{c}(\overline{\a})$.
It suffices to show that if the pair $(R,\overline{a}^t)$ is F-pure for some $t>0$, then the pair $(R,\a^{t-\epsilon})$ is also F-pure for every $t>\epsilon>0$. 
By definition, there exists an integer $s>0$ such that $\overline{\a}^{n+s}=\a^n \overline{\a}^s$ for any integer $n \ge 0$.
Take sufficiently large $q=p^e \gg 0$ such that $\lfloor t(q-1) \rfloor \ge \lfloor (t-\epsilon)(q-1) \rfloor +s$.
If $(R,\overline{a}^t)$ is F-pure, then there exists an element $d \in \overline{\a}^{\lfloor t(q-1) \rfloor} \subset \a^{\lfloor (t-\epsilon)(q-1) \rfloor}$ such that $d^{1/q}R \hookrightarrow R^{1/q}$ splits as an $R$-module homomorphism. 
This implies the pair $(R,\a^{t-\epsilon})$ is F-pure.
\end{proof}

\begin{rmk}
$(1)$ 
The pair $(R,\a^{\mathrm{c}(\a)})$ cannot be strongly F-regular by Proposition \ref{basic} $(4)$, and it is not necessarily even F-pure. See the below example. 

$(2)$
While the strong F-regularity of the pair $(R,\a^t)$ is equivalent to the strong F-regularity of the pair $(R,\overline{\a}^t)$ (cf. \cite[Proposition 1.11]{HY}), $(R,\a^t)$ is not necessarily F-pure if $(R,\overline{\a}^t)$ is F-pure.
For example, let $R=k[[X,Y,Z]]$ be a three-dimensional complete regular local ring of characteristic two and $\a=(X^2,Y^2,Z^2) \subset R$. 
Then $\overline{\a}=(X,Y,Z)^2$ and the pair $(R,\overline{\a}^{3/2})$ is F-pure, while the pair $(R,\a^{3/2})$ is not F-pure.
%However $\mathrm{c}(\a)=\mathrm{c}(\overline{\a})=3/2$ (see Proposition \ref{basic} (6)).
\end{rmk}

\begin{eg}\label{ex1}
\renewcommand{\labelenumi}{(\roman{enumi})}
\begin{enumerate}
\item Let $R=k[[X_1, \dots, X_d]]$ be a $d$-dimensional complete regular local ring over a field $k$ of characteristic $p>0$ and $\m=(X_1, \dots, X_d)$ denotes the maximal ideal of $R$.
Since $(X_1 \dots X_d)^{(q-1)/q} \in \m^{d(q-1)/q}$ is a part of free basis of $R^{1/q}/R$, we have $\mathrm{c}(\m)=d$.
\item Let $(R,\m)$ be as in (i), $S=R^{(r)}$ the $r$-th Veronese subring of $R$ and $\mathfrak{n}$ the maximal ideal of $S$.
Then $\mathrm{c}(S,\mathfrak{n})=\mathrm{c}(R,\mathfrak{n}R)=\mathrm{c}(R,\m^r)=d/r$.
\item Let $R=k[[X,Y]]$ be a two-dimensional complete regular local ring over a field $k$ of characteristic $p>0$ and $f=X^a+Y^b \in R$ with $a,b \ge 2$. Set $r=\lceil p/a \rceil + \lceil p/b \rceil -1$. 
Then, by Lemma \ref{Fedder},  
$$\frac{r-1}{p} \le c(f) \le \min \left\{\frac{r}{p}, \frac{1}{a}+\frac{1}{b} \right\}.$$ 
In particular, when $p \equiv 1 \mod ab$, we have $\mathrm{c}(f) = 1/a+1/b$. 
%By Lemma \ref{Fedder}, the pair $(R,(f)^s)$ is F-pure if and only if $f^{\lfloor s(q-1) \rfloor} \not\in (X^q, Y^q)$ for all large $q=p^e \gg 0$. Thus $\mathrm{c}((f))=1/a+1/b$. 
\item Let $R=\oplus_{n \ge 0}R_n$ be an F-finite strongly F-regular Gorenstein graded ring generated by $R_1$ over a field $R_0$, and denote $R_+=\oplus_{n \ge 1}R_n$. 
Then $\mathrm{c}(R_+)=-a(R)$, where $a(R)$ is the a-invariant of $R$ (see \cite{GW}). 
\end{enumerate}
\end{eg}

\begin{eg}[Du Val Singularities]\label{ex2}
Let $k$ be a field of characteristic $p>0$. 
The following is the list of F-pure thresholds $\mathrm{c}(\m)$ of the maximal ideal $\m$ of Du Val singularities $k[[X,Y,Z]]/(f)$.
$$
\begin{tabular}{c|lclcl}
 & \qquad\quad $f$ & & & $\mathrm{c}(\m)$ &  \\
\hline
$\mathrm{A}_n$ & $XY+Z^{n+1}$ & $1$ & \\
$\mathrm{D}_n$ & $X^2+Y^{n-1}+YZ^2$ &  $\frac{1}{2}$ &$(p>2)$\\
$\mathrm{E_6}$ & $X^2+Y^3+Z^4$ &  $\frac{1}{3}$ &$(p \equiv 1 \mod 6)$, & $\frac{1}{3}-\frac{1}{6p}$ & $(p \equiv 5 \mod 6)$ \\
$\mathrm{E_7}$ & $X^2+Y^3+YZ^3$ & $\frac{1}{4}$ &$(p \equiv 1 \mod 4)$, & $\frac{1}{4}-\frac{1}{4p}$ & \textup{($p\equiv 3 \mod 4$ and $p >3$)} \\
$\mathrm{E_8}$ & $X^2+Y^3+Z^5$& $\frac{1}{6}$ &$(p \equiv 1 \mod 6)$, & $\frac{1}{6}-\frac{1}{3p}$ & \textup{($p\equiv 5 \mod 6$ and $p > 5$)}
\end{tabular}
$$
\end{eg}
\begin{proof}
We will show the $\mathrm{E_8}$ case only, because the other cases also follow
from a similar argument. By Lemma \ref{Fedder}, the pair $(k[[X,Y,Z]]/(f),
\m^t)$ is F-pure if and only if $f^{q-1}(X,Y,Z)^{\lfloor t(q-1)
\rfloor} \not\subset (X^q, Y^q, Z^q)$ for all large $q=p^e \gg 0$. 
%where $x,y,z$ are the images of $X,Y,Z$ in $k[[X,Y,Z]]/(f)$.
If $p \equiv 1 \mod 6$, then the term $X^{2 \cdot (q-1)/2}Y^{3 \cdot (q-1)/3}Z^{5 \cdot (q-1)/6}$ appears in the expansion of
$f^{q-1}$.
Hence the F-pure threshold $\mathrm{c}(\m)$ is the supremum of all real
numbers $t >0$ such that
$$X^{q-1}Y^{q-1}Z^{\frac{5}{6}(q-1)}(X,Y,Z)^{\lfloor t(q-1)\rfloor}\not\subset (X^q, Y^q, Z^q).$$
Thus $\mathrm{c}(\m)=1-5/6=1/6$. \par
Now assume that $p \equiv 5 \mod 6$ and $p >5$, and put $q=p^e$.
By the expansion $q-1 = (p-1)p^{e-1} + (p-1)p^{e-2} + \dots + (p-1)$,
we have
$$f^{q-1} = (X^{2p^{e-1}}+Y^{3p^{e-1}}+Z^{5p^{e-1}})^{p-1}
(X^{2p^{e-2}}+Y^{3p^{e-2}}+Z^{5p^{e-2}})^{p-1} \cdots (X^{2}+Y^{3}+Z^{5})^{p-1}.$$
Note that in the expansion above, every binomial coefficient
is nonzero.
We also note that the monomial of the lowest degree in
$(X^{2p^{e-1}}+Y^{3p^{e-1}}+Z^{5p^{e-1}})^{p-1}$
which is not in $(X^q,Y^q,Z^q)$ is
$X^{2p^{e-1} \cdot (p-1)/2}Y^{3p^{e-1} \cdot (p-2)/3}Z^{5p^{e-1} \cdot (p+1)/6}$. 
So the monomial of the lowest degree in the
expansion above which is not in $(X^q,Y^q,Z^q)$ is the term
%$$(X^{(p-1) p^{e-1}}Y^{(p-2) p^{e-1}} Z^{\frac{5(p+1)p^{e-1}}{6}})(X^{(p-1)p^{e-2}}Y^{\frac{3(p-1)p^{e-2}}{2}})(X^{(p-1)p^{e-3}}Y^{3(p-1)p^{e-3}/2})\cdots,$$
$$(X^{(p-1) p^{e-1}}Y^{(p-2)p^{e-1}} Z^{\frac{5(p+1)}{6}p^{e-1}})(X^{(p-1)p^{e-2}}Y^{\frac{3(p-1)}{2}p^{e-2}}) \cdots (X^{p-1}Y^{\frac{3(p-1)}{2}}),$$
producing $X^{q-1}Y^{q-1-b_q}Z^{q-1-c_q}$, where $b_q= (p^{e-1}+1)/2$ and
$c_q= q/6 - 5p^{e-1}/6 -1$. Thus the monomial $M_q$ of the highest degree
with $M_q f^{q-1} \not\in (X^q,Y^q,Z^q)$ is $Y^{b_q} Z^{c_q}$ and we have
$$\mathrm{c}(\m) = \lim_{q\to \infty} \frac{b_q+c_q}{q-1}=\lim_{q\to \infty}\frac{q/6-p^{e-1}/3 - 1/2}{q-1} = \frac{1}{6} -\frac{1}{3p}.$$
\end{proof}
%We only show the $\mathrm{E_6}$ case, because the other cases also follow from a similar argument. By Lemma \ref{Fedder}, the pair $(k[[X,Y,Z]]/(f), (x,y,z)^t)$ is F-pure if and only if $f^{q-1}(x,y,z)^{\lfloor t(q-1) \rfloor} \not\subset (x^q, y^q, z^q)$ for all large $q=p^e \gg 0$, where $x,y,z$ are the images of $X,Y,Z$ in $k[[X,Y,Z]]/(f)$. 
%If $p \equiv 1 \mod 6$, then the term $x^{2 \cdot (q-1)/2}y^{3 \cdot (q-1)/3}z^{4\{(q-1)-(q-1)/2-(q-1)/3\}}$ appears in the expansion of $f^{q-1}$. 
%Hence the F-pure threshold $\mathrm{c}(\m)$ is the supremum of all real numbers $t >0$ such that 
%$$x^{q-1}y^{q-1}z^{\frac{2}{3}(q-1)}(x,y,z)^{\lfloor t(q-1) \rfloor}\not\subset (x^q, y^q, z^q).$$ 
%Thus $\mathrm{c}(\m)=1-2/3=1/3$.
%\end{proof}

\begin{prop}\label{key lemma}
Let $(R,\m)$ be an F-finite F-pure local ring of characteristic $p>0$ with infinite residue field $R/\m$ and $I \subsetneq R$ a proper ideal of positive height.
\begin{enumerate}
\item $\mathrm{c}(I) \le \hight I$.
\item 
Let $J \subset R$ be any minimal reduction of $I$ and denote by $\lambda>0$ the analytic spread of $I$, that is, the minimal number of generators of $J$. 
If $\mathrm{c}(I) > \lambda -r$ for an integer $r \ge 1$, then we have an inclusion $I^r \subset J$.
\end{enumerate}
\end{prop}

\begin{proof}
(1) 
We will show that the pair $(R, I^t)$ cannot be F-pure for any real number  $t > \hight I$. 
Fix a real number $t >\hight I$. 
Considering the localization at a minimal prime ideal of $I$, we may assume that $I$ is $\m$-primary. 
Since the residue field $R/\m$ is infinite, $I$ has a reduction ideal $J$ generated by at most $\hight I$ elements (see \cite{NR}). 
By definition, there exists an integer $s>0$ such that $I^{n+s}=I^sJ^n$ for any integer $n \ge 0$.
When we take sufficiently large $q=p^e \gg 0$ so that $\lfloor t(q-1) \rfloor \ge \hight I \cdot q+s$, there exists no element $d \in I^{\lfloor t(q-1) \rfloor}$ such that $d^{1/q}R \hookrightarrow R^{1/q}$ splits as an $R$-linear map, because $I^{\lfloor t(q-1) \rfloor} \subset J^{\hight I \cdot q} \subset J^{[q]}$. 
Thus the pair $(R,I^t)$ is not F-pure.

(2) 
By definition, there exists an integer $s>0$ such that $I^{n+s}=I^sJ^n$ for every integer $n \ge 0$.
Take sufficiently large $q=p^e \gg 0$ such as $(q-1)\mathrm{c}(I) > (\lambda-r)q+s$.
By the definition of the F-pure threshold, there exist an element $d \in I^{(\lambda-r)q+s}$ and an $R$-module homomorphism $\phi:R^{1/q} \to R$ sending $d^{1/q}$ to $1$. 
Now fix an arbitrary element $x \in I^r$.
Then $dx^q \in I^{\lambda q+s} \subset J^{\lambda q} \subset J^{[q]}$, because $J$ is generated by at most $\lambda$ elements.
Therefore, we have $x=\phi(d^{1/q}x) \in \phi(JR^{1/q})=J$.
\end{proof}

Now we investigate properties of local rings $(R,\m)$ with big $\mathrm{c}(\m)$. 
\begin{thm}\label{big value}
Let $(R,\m)$ be a $d$-dimensional F-finite F-pure normal local ring of characteristic $p>0$ with infinite residue field $k=R/\m$. 
\begin{enumerate}
\item 
$R$ is regular if and only if $\mathrm{c}(\m) >d-1$. 
In this case, $\mathrm{c}(\m)=d$. 
\item If $R$ is Cohen-Macaulay $\Q$-Gorenstein of index $r \ge 2$, then $\mathrm{c}(\m) \le d-1-1/r$.
\item Assume that $R$ is Cohen-Macaulay $\Q$-Gorenstein. Then $\mathrm{c}(\m)=d-1$ if and only if there exists a regular sequence $x_1, \dots, x_{d-1} \in \m$ such that $$\widehat{R}/(x_1, \dots, x_{d-1})\widehat{R} \cong k[[X,Y]]/(XY).$$
\end{enumerate}
\end{thm}
\begin{proof}
(1) 
Suppose $\mathrm{c}(\m) >d-1$ and let $J \subset R$ be a minimal reduction of $\m$.
By Proposition \ref{key lemma} (2), we have $\m=J$.
Since $J$ is generated by at most $d$ elements, $R$ is regular. 
Conversely, we assume that $R$ is regular and take a regular system of parameters $x_1, \dots, x_d$.
Then for every $q=p^e>0$, $(x_1 \cdots x_d)^{(q-1)/q} \in \m^{d(q-1)/q}$ is a part of a free basis of $R^{1/q}/R$, hence $\mathrm{c}(\m)=d$ by Proposition \ref{key lemma} (1).

(2) 
Let $J \subset R$ be a minimal reduction of $\m$. 
Then $J$ is generated by a system of parameters and there exists an integer $s>0$ such that $\m^{n+s}=J^n\m^s$ for every integer $n \ge 0$. 
Let $\a \subset R$ be a divisorial ideal isomorphic to the canonical module $\omega_R$. 
For a generator $z$ of the socle of $\a/J\a$, there exist $y \in \m \setminus J$ and $a \in \a$ such that $z=ay \mod J\a$, because $R/J$ is the Matlis dual of $\a/J\a \cong \omega_R/J\omega_R$.

Assume to the contrary that $\mathrm{c}(\m)>d-1-1/r$. 
Take sufficiently large $q=p^e \gg 0$, then by Lemma \ref{injective}, there exists an element $d \in \m^{\lceil (d-1-1/r)q \rceil+s}$ such that $dF^e:H^d_{\m}(\omega_R) \to H^d_{\m}(\omega_R^{(q)})$ is injective, where $F^e:H^d_{\m}(\omega_R) \to H^d_{\m}(\omega_R^{(q)})$ is the $e$-times iterated Frobenius map induced on $H^d_{\m}(\omega_R)$. 
By the following commutative diagram, $dz^q=d(ay)^q \mod J^{[q]}\a^{(q)} \ne 0$ in $\a^{(q)}/J^{[q]}\a^{(q)}$.
$$\xymatrix{
\a/J\a \ar@{^{(}->}[r] \ar[d]_{dF^e} & H^d_{\m}(\omega_R) \ar[d]_{dF^e} \\
\a^{(q)}/J^{[q]}\a^{(q)} \ar@{^{(}->}[r] & H^d_{\m}(\omega_R^{(q)}).\\
}$$

On the other hand $\a^r \subset \m \a^{(r)}$, because $R$ is $\Q$-Gorenstein of index $r \ge 2$. 
Hence 
$$d(ay)^q \in \m^{\lceil (d-1-1/r)q \rceil +s} \cdot \m^{\lfloor q/r \rfloor}\a^{(q)} \cdot \m^q =\m^{dq+s}\a^{(q)} \subset J^{[q]}\a^{(q)}.$$
This is a contradiction.

(3) 
First suppose that there exists a regular sequence $x_1, \dots, x_{d-1} \in \m$ such that $\widehat{R}/(x_1, \dots, x_{d-1})\widehat{R} \cong k[[X,Y]]/(XY)$. 
Then $\widehat{R}/(x_1, \dots, x_{d-2})\widehat{R}$ is the Du Val singularity of type $\mathrm{A}_n$, and by Example \ref{ex2} and the repeated applications of Proposition \ref{restriction}, $\mathrm{c}(R,\m) =\mathrm{c}(\widehat{R}, \m\widehat{R})\ge d-1$. 
Since $R$ is not regular, by (1), we have $\mathrm{c}(R,\m)=d-1$.

Now we will prove the converse implication.
Suppose that $\mathrm{c}(\m)=d-1$, then $R$ is a non-regular Gorenstein ring  by $(1)$ and $(2)$. 
Let $J \subset R$ be a minimal reduction of $\m$, then $\m^2 \subset J$ by Proposition \ref{key lemma} $(2)$. Since $J$ is generated by a system of parameters, $R/J$ is also Gorenstein. It follows from $(\m/J)^2=0$ that the embedding dimension of $R/J$ is one. 
Thus the $\m$-adic completion $\widehat{R}$ is a hypersurface singularity of degree two.  
If there exists no regular sequence $x_1, \dots, x_{d-1} \in \m$ satisfying $\widehat{R}/(x_1, \dots, x_{d-1})\widehat{R} \cong k[[X,Y]]/(XY)$, then $\widehat{R} \cong k[[X,Y_1, \dots, Y_d]]/(X^2+g(Y_1, \dots, Y_{d}))$, where $g(Y_1, \dots, Y_d)$ is an element in $k[[Y_1, \dots, Y_{d}]]$ of degree $\ge 3$. 
We write by $x, y_1, \dots, y_{d}$ the images of $X, Y_1, \dots, Y_{d}$ in $\widehat{R}$. 
We can assume that $J\widehat{R}=(y_1, \dots, y_d) \subset \widehat{R}$, then the image of $x$ is a generator of the socle of $\widehat{R}/J\widehat{R}$. 
Since $\mathrm{c}(R,\m)=d-1>d-3/2$, by an argument similar to the proof of (2), for all large $q=p^e \gg 0$, there exists an element $d \in (x,y_1, \dots, y_{d})^{\lfloor (d-3/2)(q-1) \rfloor+2}$ such that $dx^q \not\in (J\widehat{R})^{[q]}$.
Since $x^2=g(y_1, \dots, y_{d})$ in $\widehat{R}$, we have $x^q \in (J\widehat{R})^{3\lfloor q/2 \rfloor}$ and $d \in (J\widehat{R})^{\lfloor (d-3/2)(q-1) \rfloor +1}$. 
Therefore, $dx^q \in (J\widehat{R})^{d(q-1)+1} \subset (J\widehat{R})^{[q]}$. 
This is a contradiction, hence we obtain the assertion.
\end{proof}

\begin{rmk}\label{d-1}
When $R$ is not $\Q$-Gorenstein, Theorem \ref{big value} (3) does not hold true. 
For example, let $R=k[[X_{ij} \mid 1 \le i \le d-1, j=1,2]]/(X_{k1}X_{l2}-X_{k2}X_{l1} \mid 1 \le k \le l \le d-1)$ be the completion of the Segre product of two-dimensional and $d-1$-dimensional polynomial rings over a field $k$ of characteristic $p>0$. 
%the Segre product of a $(d-1)$-dimensional and a two-dimensional formal power s%eries ring 
Then $R$ is not $\Q$-Gorenstein but $\mathrm{c}(R,\m)=d-1$.
\end{rmk}

We expect the converse of Remark \ref{d-1}.

\begin{conj}
Let $(R,\m)$ be a $d$-dimensional F-finite F-pure non-$\Q$-Gorenstein normal local ring of characteristic $p>0$. 
If $\mathrm{c}(\m)=d-1$, then the associated graded ring $\mathrm{gr}_{\m}(R)$ of $R$ with respect to $\m$ is isomorphic to the ring $k[X_{ij} \mid 1 \le i \le d-1, j=1,2]/(X_{k1}X_{l2}-X_{k2}X_{l1} 
\mid 1 \le k \le l \le d-1)$.
\end{conj}

When $R$ is $\Q$-Gorenstein Cohen-Macaulay, we have a more powerful proposition of Proposition \ref{key lemma} $(2)$, which has an interesting application to algebraic geometry (see Proposition \ref{kakimi}). Before we state the result, let us introduce some notation. 
For any ideal $\b \subset R$ (resp. any element $x \in R$), we define the {\it vanishing order} $\mathrm{ord}_{\a}(\b)$ (resp. $\mathrm{ord}_{\a}(x)$) of $\b$ 
(resp. $x$) with respect to an ideal $\a$ to be the supremum of 
all integers $k$ such that $\b \subset \a^k$ (resp. $x\in \a^k$). 
\begin{prop}\label{reduction}
Let $(R,\m)$ be an F-finite F-pure normal local ring of characteristic $p>0$ with infinite residue field $R/\m$, $J \subset R$ a minimal reduction of $\m$ and $\a \subset R$ a divisorial ideal isomorphic to the canonical module $\omega_R$. 
Assume that $R$ is Cohen-Macaulay $\Q$-Gorenstein of index $r \ge 2$ and the symbolic $r$-th power $\a^{(r)}$ of $\a$ is generated by an element $x \in \m$. 
If $\mathrm{c}(\m)>d-k-1+\mathrm{ord}_{\m}(x)/r$ for an integer $k \ge 1$, then $\m^k \subset J$. 
\end{prop}
\begin{proof}
The proof is similar to that in Proposition \ref{big value} $(2)$. 
Assume to the contrary that $\m^k \not\subset J$. 
Since $\a/J\a$ is the Matlis dual of $R/J$, there exist $y \in \m^k \setminus J$ and $a \in \a$ such that $ay \mod J\a \in \a/J\a$ is a generator of the socle of $\a/J\a$. 
Since $J$ is a minimal reduction of $\m$, there exists an integer $s>0$ such that $\m^{n+s}=J^n\m^s$ for every integer $n \ge 0$. 
Let $l=r-\mathrm{ord}_{\m}(\a^{(r)})>0$. 
Since $\mathrm{c}(\m)>d-k-l/r$, when we take sufficiently large $q=p^e \gg 0$, there exists $d \in \m^{\lceil (d-k-l/r)q \rceil +s}$ such that $d^{1/q}R \hookrightarrow R^{1/q}$ splits as an $R$-linear map. 
Hence, by an argument similar to the proof of Proposition \ref{big value} $(2)$, $d(ay)^q$ does not belong to $J^{[q]}\a^{(q)}$. 
However, since there exists some element $b \in \m^l$ such that $a^r=bx$, 
$$d(ay)^q \in \m^{\lceil (d-k-l/r)q \rceil +s} \cdot \m^{l\lfloor q/r \rfloor}\a^{(q)} \cdot \m^{kq} =\m^{dq+s}\a^{(q)} \subset J^{[q]}\a^{(q)}.$$
This is a contradiction, hence $\m^k$ is contained in $J$.
\end{proof}

\section{F-pure threshold in characteristic zero}
The F-pure threshold is also defined for an ideal of a ring of characteristic zero. 

\begin{dfn}\label{char 0}
Let $R$ be a reduced algebra essentially of finite type over a field $k$ of characteristic zero and $\a \subset R$ an ideal such that $\a \cap R^{\circ} \neq \emptyset$. 
Suppose that $R$ is of dense F-pure type. Then 
the {\it F-pure threshold} $\mathrm{c}(\a)$ of $\a$ is defined to be 
$$\mathrm{c}(\a)=\mathrm{c}(R,\a)=\sup \{s \in \R_{\ge 0} \mid \textup{the pair $(R,\a^s)$ is of dense F-pure type} \}.$$
\end{dfn}

The F-pure threshold $\mathrm{c}(\a)$ in characteristic zero corresponds to the lc threshold $\mathrm{lc}(\a)$ under mild conditions.  
In particular, when $R$ is a log terminal singularity of characteristic zero, we have $\mathrm{c}(\a)=\mathrm{lc}(\a)$. 
\begin{prop}\label{lc thresholds}
Let $R$ be a $\Q$-Gorenstein normal local ring of essentially of finite type over a field of characteristic zero and $\a \subset R$ a nonzero ideal. 
Suppose that $R$ is of dense F-pure type. 

\begin{enumerate}
\item $\mathrm{c}(\a) \le \mathrm{lc}(\a)$. 
\item Let $I \subset R$ be the radical ideal defining the non-log-terminal locus of $\Spec R$. 
If $I \supset \a$, then $\mathrm{c}(\a)=\mathrm{lc}(\a)$. 
\end{enumerate}
\end{prop}
\begin{proof}
$(1)$ By Proposition \ref{klt} $(1)$, it is obvious. 

$(2)$ Since $I \supset \a$, by Remark \ref{multiplier ideals} $(4)$, 
$$\mathrm{lc}(\a)=\sup\{s \in \R_{\ge 0} \mid \textup{the pair $(\Spec R,sV(\a))$ is klt} \}. $$
Since klt pairs are equivalent to $\Q$-Gorenstein pairs of strongly F-regular type by Proposition \ref{klt} $(2)$, the assertion follows from Proposition \ref{basic} $(5)$.
\end{proof}

The following proposition characterizes local rings $(R,\m)$ satisfying $\mathcal{J}(\a) \supsetneq \a$ for every $\m$-primary ideal $\a \subset R$ (cf. \cite[Corollary 2.3]{HT}). 
\begin{prop}
Let $(R,\m)$ be a $\Q$-Gorenstein normal local ring essentially of finite type over a field of characteristic zero and $\a \subset R$ an $\m$-primary ideal. 
Assume that $R$ is of dense F-pure type and the non-log-terminal locus of $\Spec R$ is isolated or empty.
Then $\mathrm{c}(\a) >1$ if and only if for every $\m$-primary ideal $\b \subset R$, we have $J(\b) \supseteq \b:\a$.
In particular, $\mathrm{c}(\m) >1$ if and only if for every $\m$-primary ideal $\b \subset R$, we have a strict containment $\J(\a) \supsetneq \a$. 
\end{prop}
\begin{proof}
First suppose that $J(\b) \supseteq \b:\a$ for every $\m$-primary ideal $\b \subset R$. 
Considering the case where $\a=\b$, we have $\J(\a)=R$. 
That is, the pair $(\Spec R, V(\a))$ is klt.
Hence $\mathrm{c}(\a)>1$ by Proposition \ref{lc thresholds} (2).

Conversely assume $\mathrm{c}(\a)>1$. 
By Proposition \ref{lc thresholds} (2), this implies that $\J(\a)=R$. 
By the definition of multiplier ideals, there exists a log resolution $f:X \to \Spec R$ such that the ideal sheaf $\a \O_X=\O_X(-Z_a)$ is invertible and $\lceil K_{X/\Spec R} \rceil \ge Z_a$. 
We may assume that $\b \O_X=\O_X(-Z_b)$ is also invertible and $\b$ is an integrally closed ideal, that is, $\b=H^0(X,\O_X(-Z_b))$. 
Fix an arbitrary element $x \in \b:\a$. 
By definition, for any element $a \in \a$, we have $ax \in H^0(X,\O_X(-Z_b))$, namely, $\Div_X(a)+\Div_X(x) \ge Z_b$. 
Note that by taking the element $a$ sufficiently general, $\Div_X(a)=Z_a+f^{-1}_*\Div_R(a)$. 
Hence, looking at the exceptional part, we have 
\begin{align*}
\Div_X(x)+\lceil K_{X/\Spec R} \rceil-Z_b &\ge \lceil K_{X/\Spec R}\rceil -Z_a\\
&\ge 0.
\end{align*}
Thus $x \in \J(\b)$.
\end{proof}

\begin{eg}
Let $(R,\m)$ be a $d$-dimensional regular local ring essentially of finite type over a field of characteristic zero. 
Since $\mathrm{c}(\m^{d-1})=d/(d-1) >1$, we have $\J(\a) \supset (\a:\m^{d-1})$ for every $\m$-primary ideal $\a \subset R$.
\end{eg}

In the latter half of this section, we investigate properties of three-dimensional terminal singularities in terms of the F-pure threshold. 

\begin{prop}\label{Gorenstein terminal}
Let $(R,\m)$ be a three-dimensional Gorenstein normal local ring essentially of finite type over a field of characteristic zero and $I \subset R$ the radical ideal defining the singular locus of $\Spec R$. 
Suppose that $R$ is of dense F-pure type. 
Then $\Spec R$ has at worst terminal singularities if and only if $\mathrm{c}(I) >1$. 
\end{prop}

\begin{proof}
First assume that $R$ is a terminal singularity. 
By \cite{Re1}, a three-dimensional normal singularity is terminal of index one if and only if it is an isolated cDV singularity. 
Thus $I=\m$ and $R/xR$ is a Du Val singularity for some $x \in \m$. 
Denote $S=R/xR$, then by Example \ref{ex2}, $\mathrm{c}(S,IS)>0$. 
Therefore, by Proposition \ref{restriction}, we obtain the inequality $\mathrm{c}(R,I) \ge \mathrm{c}(S, IS)+1>1$. 

Conversely, we suppose that $\mathrm{c}(I)>1$. 
Let $f:\widetilde{X} \to X=\Spec R$ be a log resolution of $I$ so that $I\O_{\widetilde{X}}=\O_{\widetilde{X}}(-F)$ is invertible. 
Since $\mathrm{lc}(I)=\mathrm{c}(I)>1$ by Proposition \ref{lc thresholds} $(2)$, we have $\mathcal{J}(I)=R$, that is, $\lceil K_{\widetilde{X}/X} \rceil \ge F$. 
This implies that the coefficient of $K_{\widetilde{X}/X}$ in each irreducible component is greater than zero, because $F$ is supported on the exceptional locus of $f$ by the definition of $I$. 
Therefore, $\Spec R$ has at worst terminal singularities. 
\end{proof}

\begin{rmk}
The proof above shows that if $\mathrm{c}(I)>1$, then $R$ is a terminal singularity 
in any dimension. 
In \cite{Wa}, the second-named author called strongly F-regular rings with 
 $\mathrm{c}(I) > 1$ F-terminal rings, where $I$ is the defining ideal of the 
singular locus of $R$. 
This property is preserved under localization and F-terminal rings are 
regular in codimension two. 
In \cite{Wa}, it is asserted that terminal singularities in dimension three 
are of F-terminal type by an erroneous proof. 
But as is shown in Example \ref{2/3}, there are non-Gorenstein terminal 
singularities in dimension three, with $c(\m) \le 1$. 
By this reason, we withheld to use the term ``F-terminal" in this 
paper.
\end{rmk}

If the dimension is greater than three, the condition that $\mathrm{c}(\m) >1$ is stronger than being a terminal singularity.
\begin{eg}[\textup{\cite{R}}]
Let $R=\C[[X,Y,Z,U,V]]/(X^3+Y^4+Z^4+U^6+V^6)$ and denote by $\m$ the maximal ideal of $R$. 
Then $R$ is an isolated Gorenstein terminal singularity but $\mathrm{c}(\m)=1$. 
\end{eg}

%Also, if $R$ is a three-dimensional $\Q$-Gorenstein normal local ring of index $\ge 2$, we cannot give a simple characterization as Proposition \ref{Gorenstein terminal}. 
%We need the following proposition to compute the F-pure threshold $\mathrm{c}(\m)$ of the maximal ideal $\m$ of a non-Gorenstein local ring $(R,\m)$. 
Now, we will compute a lower bound for $\mathrm{c}(\m)$ of  three-dimensional 
$\Q$-Gorenstein terminal singularities  of index $\ge 2$. 

\begin{prop}\label{lemma for terminal}
Let $(R,\m)$ be a $d$-dimensional $\Q$-Gorenstein Cohen-Macaulay normal local ring of characteristic $p>0$ with Gorenstein index $r>0$. 
Let $\a \subset R$ be a divisorial ideal isomorphic to the canonical module $\omega_R$. Then the symbolic $r$-th power $\a^{(r)}$ of $\a$ is a principal ideal generated by an element $w \in \m$. 
If $R/\a$ is Gorenstein strongly F-regular, then 
$$\mathrm{c}(R,\m) \ge \mathrm{c}(R/\a, \m/\a)+\frac{\mathrm{ord}_{\m}(w)}{r}.$$ 
\end{prop}

\begin{proof}
Take a system of parameters $x_1, x_2, \dots, x_d$ of $R$ such that $x_1 \in \a$ and the image of $x_2, \dots, x_d$ is a system of parameters of $R/\a$. 
Set $\mbox{\boldmath $x$}=\{x_2, \dots, x_d\}$, and for an integer $n>0$, write $\mbox{\boldmath $x$}^n=\{x_2^n, \dots, x_d^n\}$. 
Suppose that the pair $(R/\a, (\m/\a)^c)$ is strongly F-regular for some $c>0$, and then we will prove that $(R,\m^{c+\mathrm{ord}_{\m}(w)/r})$ is an F-pure pair. 
%Let $z \in R$ be an element whose image in $R/(\a+(x_2, \dots, x_d))$ generates the socle of $R/(\a+(x_2, \dots, x_d))$. 
%Since the pair $(R/\a, (\m/\a)^c)$ is strongly F-regular, by a similar argument in Proposition \ref{big value} (2), for all large $q=p^e>0$, there exists an element $d \in \m^{\lceil cq \rceil}$ such that $dwz^q \not\in I+(x_2^q, \dots, x_d^q)$. 
By Lemma \ref{injective}, for all large $q=p^e \gg 0$, there exists an element $d \in \m^{\lceil cq \rceil}$ such that $dwF^e_{R/\a}:H^{d-1}_{\m/\a}(R/\a) \to H^{d-1}_{\m/\a}(R/\a)$ is injective, where $F^e_{R/\a}:H^{d-1}_{\m/\a}(R/\a) \to H^{d-1}_{\m/\a}(R/\a)$ is the $e$-times iterated Frobenius map induced on $H^{d-1}_{\m/\a}(R/\a)$. 
Take an element $z \in R$ whose image in $R/(\a+(\mbox{\boldmath $x$}))$ generates the socle of $H^{d-1}_{\m/\a}(R/\a)$. 
Since $H^{d-1}_{\m/\a}(R/\a) = \varinjlim R/(\a+(\mbox{\boldmath $x$}^s))$ (the direct limit map $R/(\a+(\mbox{\boldmath $x$}^s)) \to R/(\a+(\mbox{\boldmath $x$}^{s+1}))$ is the multiplication by $x_2 \cdots x_d$), $dwz^q$ does not belong to $\a+(\mbox{\boldmath $x$}^q)$ for all large $q=p^e \gg 0$. 

By Lemma \ref{injective} again, the pair $(R,\m^{c+\mathrm{ord}_{\m}(w)/r})$ is F-pure if and only if for all large $q=p^e \gg 0$, there exists an element $d' \in \m^{\lfloor (c+\mathrm{ord}_{\m}(w)/r)(q-1) \rfloor}$ such that $d'F^e_R:H^{d}_{\m}(\omega_R) \to H^{d}_{\m}(\omega_R^{(q)})$ is injective, where $F^e_R:H^{d}_{\m}(\omega_R) \to H^{d}_{\m}(\omega_R^{(q)})$ is the $e$-times iterated Frobenius map induced on $H^{d}_{\m}(\omega_R)$. 
Note that $H^{d}_{\m}(\omega_R) \cong \varinjlim \a/(x_1^s, \mbox{\boldmath $x$}^s)\a$, $H^{d}_{\m}(\omega_R^{(q)}) \cong \varinjlim \a^{(q)}/(x_1^{sq}, \mbox{\boldmath $x$}^{sq})\a^{(q)}$ and the image of $x_1z$ in $\a/(x_1, \mbox{\boldmath $x$})\a$ is a generator of the socle of $H^{d}_{\m}(\omega_R)$.
Write $q=kr+i$ for integers $k$ and $i$ with $0 \le i \le r-1$. 
Since $dw^{k+1} \in \m^{\lceil cq \rceil + \mathrm{ord}_{\m}(w)(k+1)} \subseteq \m^{\lfloor (c+\mathrm{ord}_{\m}(w)/r)(q-1) \rfloor}$, it is sufficient to show that $dw^{k+1}F^e_R$ is injective, that is, $dw^{k+1}(x_1z)^q \not\in (x_1^q, \mbox{\boldmath $x$}^q)\a^{(q)}$ for all large $q=p^e \gg 0$.
If $dw^{k+1}(x_1z)^q \in (x_1^q, \mbox{\boldmath $x$}^q)\a^{(q)}=w^k(x_1^q, \mbox{\boldmath $x$}^q)\a^{(i)}$, then there exists an element $a \in \a^{(i)}$ such that $(dwz^q-a)x_1^q \in (\mbox{\boldmath $x$}^q)\a^{(i)}$, hence $dwz^q \in \a+(\mbox{\boldmath $x$}^q)$.
This is a contradiction.
\end{proof}

Let $(0 \in X)$ be a three-dimensional terminal singularity of index $r \ge 2$ over the complex number field $\C$.  
Using the classification of three-dimensional terminal singularities due to Mori \cite{M} and Reid \cite{R}, we compute the F-pure threshold $\mathrm{c}(\O_{X,0},\m_{X,0})$.  
Let $\pi:(0 \in \widetilde{X}) \to (0 \in X)$ be an index one cover. 
Then the group $\mu_r$ of $r$-th roots of unity acts on $\widetilde{X}$ and a general element $H \in |-K_X|$ which contains $0$ is Du Val. 
Applying Proposition \ref{lemma for terminal} to the list of $\tilde{H}=\pi^{*}(H)$, $H$ and the action of $\mu_r$ on $\C^4$, we give a lower bound for the F-pure threshold $\mathrm{c}(\O_{X,0},\m_{X,0})$.  
$$
\begin{tabular}{l|llll}
 & $\quad \tilde{H} \to H$ & $r$ & the action of $\mu_r$ & \textup{$\mathrm{c}(\O_{X,0}, \m_{X,0})$} \\
\hline
$\mathrm{cA/}r$ & $\mathrm{A}_{n-1} \to \mathrm{A}_{rn-1}$ & $r$ & $1/r(a,-a,1,0;0)$ & $\ge 1+1/r$ \\ 
$\mathrm{cAx/4}$ & $\mathrm{A}_{2n-2} \to \mathrm{D}_{2n+1}$ & $4$ & $1/4(1,1,3,2;2)$ & $\ge 3/4$ \\
$\mathrm{cAx/2}$ & $\mathrm{A}_{2n-1} \to \mathrm{D}_{n+2}$ & $2$ & $1/2(0,1,1,1;0)$ & $\ge 1$ \\
$\mathrm{cD/3}$ & $\mathrm{D}_{4} \to \mathrm{E}_6$ & $3$ & $1/3(0,2,1,1;0)$ & $\ge 2/3$ \\
$\mathrm{cD/2}$ & $\mathrm{D}_{n+1} \to \mathrm{D}_{2n}$ & $2$ & $1/2(1,0,1,1;0)$ & $\ge 1$ \\
$\mathrm{cE/2}$ & $\mathrm{E}_6 \to \mathrm{E}_7$ & $2$ & $1/2(1,0,1,1;0)$ & $\ge 3/4$ \\
\end{tabular}
$$
In the above list, $1/r(a_1,\dots,a_4;b)$ means that the generator $\xi$ of $\mu_r$ acts on the coordinates $x_1, \dots, x_4$ and on the equation $f$ of $\widetilde{X}$ as 
$$(x_1, \dots, x_4) \mapsto (\xi^{a_1}x_1, \dots, \xi^{a_4}x_4;\xi^bf).$$ 
By the above list, we know $\mathrm{c}(\O_{X,0}, \m_{X,0}) \ge 2/3$. 
Indeed, we can construct some example with $\mathrm{c}(\O_{X,0}, \m_{X,0})=2/3$. 
\begin{eg}\label{2/3}
Let $S=\C[[X,Y,Z,T]]/(X^2+Y^3+Z^3+YT^4+T^6)$ and let $\mu_3$ act on $S$ as $1/3(0,2,1,1;0)$.  
Then the invariant subring 
$$R=(\C[[X,Y,Z,T]]/(X^2+Y^3+Z^3+YT^4+T^6))^{\mu_3} \subset S$$ 
is a three-dimensional terminal singularity of type $\mathrm{cD/3}$ (see the above list). 
We will prove that $\mathrm{c}(R, \m_R)=2/3$, where $\m_R$ denotes the maximal ideal of $R$. 

Considering reduction to characteristic $p \gg 0$, we may assume that $S$ is defined over a field of characteristic $p>0$. 
Denote by $x,y,z,t$ the images of $X,Y,Z,T$ in $S$. 
Then 
$$\a=(tS)^{\mu_3}=(yt,z^2t, zt^2, t^3)R \subset R$$ 
is isomorphic to the canonical module $\omega_R$ and $R/\a$ is a Du Val singularity of type $\mathrm{E_6}$. 
Furthermore, the symbolic cube $\a^{(3)}$ of $\a$ is a principal ideal generated by $t^3$. 
Denote $J=(x, z^3, t^3) \subset R$. 
Then the image of $w=y^2z^2t^3$ is a generator of the socle of $\a/\a J$ and the image of $w_e=x^{q-1}y^2z^{3q-1}t^{4q-1}$ is a generator of the socle of $\a^{(q)}/\a^{(q)}J^{[q]}$ for every $q=p^e>0$. 
Suppose that the weights for $x, y, z, t$ are $3,2,2,1$ respectively, then the weighted degree of $w$ is $11$ and that of $w_e$ is $13q-2$. 
If the pair $(R,\m_R^t)$ is strongly F-regular for a real number $t>0$, then by Lemma \ref{injective} $(3)$, there exists an element $a \in \m_R^{\lceil tq \rceil}$ of weighted degree $d$ such $11q+d \le 13q-2$. 
Since $\m_R$ is generated by elements of weighted degree $\ge 3$, we have $t < 2/3$, hence $\mathrm{c}(R,\m_R)=2/3$.
\end{eg}

Kakimi \cite{Ka} studied the multiplicity of three-dimensional terminal singularities, using the classification due to Mori \cite{M}.
By making use of the F-pure threshold, we can give a simple unified proof of his result. 
%We believe that our proof is simpler than his, although we also base our proof on the classification theory. 
\begin{prop}\label{kakimi}
Let $(R,\m)$ be a three-dimensional $\Q$-Gorenstein normal local ring essentially of finite type over the complex number field $\C$ and $J \subset R$ a minimal reduction of $\m$. If $R$ is a terminal singularity, then $\m^2 \subset J$, particularly, 
\begin{align*}
\dim_k(\m^k/\m^{k+1})&=\mathrm{e}(R) \cdot \frac{k(k+1)}{2}+k+1,\\
\mathrm{embdim}(R)&=\mathrm{e}(R)+2,
\end{align*}
where $\mathrm{e}(R)$ denotes the multiplicity of $R$ and $\mathrm{embdim}(R)$ denotes the embedding dimension of $R$. 
\end{prop}
\begin{proof}
The Gorenstein case immediately follows from Proposition \ref{key lemma} $(2)$ and Proposition \ref{Gorenstein terminal}. 
Furthermore, when $R$ is $\Q$-Gorenstein of index $2$, applying Proposition \ref{reduction} to the above list, we have the assertion. 
Hence we may assume that $R$ is $\Q$-Gorenstein of index $r \ge 3$, that is, of type $\mathrm{cAx/4}$ or $\mathrm{cD/3}$. 
Let $S=\bigoplus_{i=0}^{r-1}\omega_R^{(i)}$ be the canonical cover of $R$, then $S$ is a hypersurface singularity $\C[X,Y,Z,T]_{(X,Y,Z,T)}/(f)$. 
We denote by $x,y,z,t$ the images of $X,Y,Z,T$ in $S$. 
Applying Proposition \ref{reduction} to the above list again, it suffices to show that for some divisorial ideal $\a \subset R$ which is isomorphic to the canonical module $\omega_R$, the symbolic $r$-th power $\a^{(r)}$ of $\a$ is generated by an element $x \notin \m^2$. 

First suppose that $R$ is of type $\mathrm{cD/3}$. 
By the above list,  $R=S^{\mu_3}$ is the localization of $\C[X,YZ,YT, Y^3, ZT^2, Z^2T, Z^3, T^3]/(f)$. 
Then the ideal $\a=(tS)^{\mu_3}=(yt, zt^2, z^2t, t^3)R \subset R$ is isomorphic to the canonical module $\omega_R$. 
Since $\a^{(3)}$ is generated by a single element, $\a^{(3)}$ has to be $(t^3)$ and $t^3$ is not contained in $\m^2$. 

Next we consider the case where $R$ is of type $\mathrm{cAx/4}$. 
Denote $\a=(xS)^{\mu_4}=(x^4, x^3y, x^2y^2, x^2t, xz, xy^3, xyt)R \subset R=S^{\mu_4}$. 
Then $\a$ is isomorphic to $\omega_R$ and $\a^{(4)}$ is generated by $x^4$. 
Since $x^4$ is not contained in $\m^2$, we obtain the assertion.  
\end{proof}

\section{Log canonical threshold vs. F-pure threshold}
In this section, we pick up some results on the lc threshold and show how we can prove these by our theory. 

%We prepare a notation. Let $R$ be a reduced ring and $\a \subsetneq R$ a proper ideal. 
%For any ideal $\b \subset R$, we define the vanishing order $\mathrm{ord}_{\a}(\b)$ of $\b$ with respect to $\a$ to be the supremum of all integers $k$ such that $\b \subset \a^k$. 
%By virtue of Proposition \ref{lc thresholds}, our results give simple proofs of those results.

First we give a fundamental bound of the F-pure threshold, which is an analogue of \cite[1.4]{DK} or \cite[Lemma 8.10]{Ko}. 
\begin{prop}\label{bound}
Let $(R,\m)$ be a $d$-dimensional F-finite regular local ring of characteristic $p>0$ and $\a \subsetneq R$ a proper ideal of positive height. 
Then 
$$\frac{1}{\mathrm{ord}_{\m}(\a)} \le \mathrm{c}(\a) \le \frac{d}{\mathrm{ord}_{\m}(\a)}. $$
\end{prop}
\begin{proof}
Denote $k=\mathrm{ord}_{\m}(\a)$. 
Since $\a^{\lceil \frac{d}{k} q \rceil} \subset (\m^k)^{\lceil \frac{d}{k} q \rceil} \subset \m^{dq} \subset \m^{[q]}$ for all $q=p^e>0$, 
by Lemma \ref{Fedder}, the pair $(R,\a^{d/k})$ is not strongly F-regular. 
Therefore, $\mathrm{c}(\a) \le d/k$. 

For the other inequality, we will show that the pair $(R,\a^{1/k-\epsilon})$ is strongly F-regular for all $1/k >\epsilon >0$. 
Take an element $x \in \a \setminus \m^{k+1}$. 
Since $\mathrm{ord}_\m(\a^{\lceil (1/k-\epsilon)q \rceil}) \le \mathrm{ord}_{\m}(x^{\lceil (1/k-\epsilon)q \rceil})<q$ for all large $q=p^e \gg 0$, by Lemma \ref{Fedder} again, the pair $(R,\a^{1/k-\epsilon})$ is strongly F-regular. 
\end{proof}

When $\a$ is generated by homogeneous elements of the same degree, we have a stronger bound for the F-pure threshold. 
The following result is a weak form of \cite[Theorem 3.4]{DEM2}.
\begin{prop}\label{bound2}
Let $S=k[X_1, \dots, X_n]$ be an $n$-dimensional polynomial ring over a perfect field of characteristic $p>0$ and $\a \subsetneq S$ a proper ideal of height $h>0$ which is generated by homogeneous elements of degree $d$. 
Then $\mathrm{c}(\a) \ge \frac{h}{d}$. 
\end{prop}
\begin{proof}
Let $f_1, \dots, f_{n-h} \in S$ be a system of linear forms whose image in $S/\a$ is a system of parameters. 
Let $\b \subset R$ be an ideal generated by $\a$ and $f_1, \dots, f_{n-h}$. 
%Then consider an $(X_1, \dots, X_n)$-primary ideal $\b=\a+(f_1, \dots, f_{n-h}) \subset S$. 
\begin{cl}
$$\mathrm{c}(\a)+n-h \ge \mathrm{c}(\b).$$
\end{cl}
\begin{proof}[Proof of Claim]
Suppose that the pair $(R,\b^t)$ is F-pure for some $t > n-h$, and then it suffices to show that the pair $(R, \a^{t-(n-h)})$ is also F-pure. 

By Remark \ref{split} and Lemma \ref{Fedder}, for all large $q=p^e \gg 0$, there exist nonnegative integers $l_0, \dots, l_{n-h}$ with $l_0+ \dots+l_{n-h}=\lfloor t(q-1) \rfloor$ such that $\a^{l_0}f_1^{l_1} \cdots f_{n-h}^{l_{n-h}} \not\subset (X_1^q, \dots, X_n^q)$. 
Since $l_i$ is clearly less than $q$ for all $i=1, \dots, n-h$, we have $\a^{\lfloor (t-(n-h))(q-1) \rfloor} \not\subset (X_1^q, \dots, X_n^q)$, that is, the pair $(R, \a^{t-(n-h)})$ is F-pure. 
\end{proof}
Since $J$ is generated by homogeneous elements of degree $d$, there exist linear forms $g_1, \dots, g_h \in S$  such that $f_1, \dots, f_{n-h}, g_1, \dots, g_h$ is a system of parameters of $S$ and that $\b$ and $(f_1, \dots, f_{n-h})+(g_1, \dots, g_h)^d$ have the same multiplicity. 
Thus $\overline{\b}=(f_1, \dots, f_{n-h})+(g_1, \dots, g_h)^d$ and the F-pure threshold of $(f_1, \dots, f_{n-h})+(g_1, \dots, g_h)^d$ is $h/d+n-h$. 
By the above claim, we obtain the inequality 
$$\mathrm{c}(\a)\ge \mathrm{c}(\b)-(n-h)=\mathrm{c}((f_1, \dots, f_{n-h})+(g_1, \dots, g_h)^d)-(n-h)=h/d.$$
\end{proof}

We show the restriction property of the F-pure threshold. 
The following is an analogue of Musta{\c{t}}{\v{a}}'s result \cite[Proposition 4.5, Corollary 4.6]{Mu} which is obtained via the theory of jet schemes. 
\begin{prop}\label{restriction}
Let $(R,\m)$ be an F-finite F-pure $\Q$-Gorenstein normal local ring of characteristic $p>0$ and $\a \subsetneq R$ a proper ideal of positive height. 
\begin{enumerate}
\item Let $x \in \m$ be a nonzero divisor of $R$ such that $R/xR$ is a normal Cohen-Macaulay F-pure ring, and denote $S=R/xR$. Then
$$\mathrm{c}(R,\a) \ge \mathrm{c}(S, \a S)+\mathrm{ord}_{\a}(x).$$

\item Assume that $R$ is a regular local ring. Let $I \subset R$ be an unmixed radical ideal such that $R/I$ is F-pure, and denote $S=R/I$. 
Then 
$$\mathrm{c}(R,\a) \ge \mathrm{c}(S, \a S)+\mathrm{ord}_{\a}(I).$$
\end{enumerate}
\end{prop}
\begin{proof}
$(1)$
Assume that the pair $(S,(\a S)^t)$ is F-pure for some $t>0$, and then it is enough to show that the pair $(R,x\a^t)$ is also F-pure. 
Let $F^e:H^{d-1}_{\m S}(\omega_S) \to H^{d-1}_{\m S}(\omega_S^{(q)})$ (resp. $F^e:H^{d}_{\m}(\omega_R) \to H^{d}_{\m}(\omega_R^{(q)})$) be the $e$-times iterated Frobenius map induced on $H^{d-1}_{\m S}(\omega_S)$ (resp. $H^d_{\m}(\omega_R))$.
Then, by Lemma \ref{injective}, for all large $q=p^e \gg 0$, there exists $d \in \a^{\lfloor t(q-1) \rfloor}$ such that $dF^e:H^{d-1}_{\m S}(\omega_S) \to H^{d-1}_{\m S}(\omega_S^{(q)})$ is injective.
By an argument similar to the proof of \cite[Theorem 4.9]{HW}, $dx^{q-1}F^e:H^{d}_{\m}(\omega_R) \to H^{d}_{\m}(\omega_R^{(q)})$ is also injective. 
By Lemma \ref{injective} again, the pair $(R,x\a^t)$ is F-pure. 

$(2)$ It suffices to show that if the pair $(S,(\a S)^t)$ is F-pure for some $t>0$, then the pair $(R,\a^tI)$ is also F-pure. 
However, it immediately follows from \cite[Theorem 3.11]{T}. 
\end{proof}

We show the summation property of the F-pure threshold with simple proof. Refer to \cite[Theorem 2.9]{DK} for the summation property of  the lc threshold. 
\begin{prop}\label{subadditivity}
Let $(R,\m)$ be an F-finite F-pure local ring of characteristic $p>0$ and $\a, \b \subseteq R$ ideals of positive height.
Then 
$$\mathrm{c}(\a)+\mathrm{c}(\b) \ge \mathrm{c}(\a+\b).$$
\end{prop}
\begin{proof}
It suffices to see that for some $t,s>0$ if neither of the pairs $(R,\a^t)$ nor $(R,\b^s)$ is strongly F-pure, then the pair $(R,(\a+\b)^{t+s})$ is also not strongly F-pure. 
Assume to the contrary that the pair $(R,(\a+\b)^{t+s})$ is strongly F-pure.
Since $(\a+\b)^{\lceil (t+s)q \rceil} \subset \a^{\lceil tq \rceil}+\b^{\lceil sq \rceil}$, the strong F-purity of the pair $(R,(\a+\b)^{t+s})$ implies that there exist $q=p^e>0$, $a \in \a^{\lceil tq \rceil}$ and $b \in \b^{\lceil sq \rceil}$ such that $(a+b)^{1/q}R \hookrightarrow R^{1/q}$ splits as an $R$-module homomorphism. 
Then either $a^{1/q}R \hookrightarrow R^{1/q}$ or $b^{1/q}R \hookrightarrow R^{1/q}$ is forced to split as an $R$-module homomorphism by Remark \ref{split}. 
That is, either the pair $(R,\a^t)$ or the pair $(R,\b^s)$ is strongly F-pure. 
\end{proof}

Finally we give an inequality between the multiplicity and the F-pure threshold of an $\m$-primary ideal, which is an analogue of \cite[Theorem 1.2, 1.4]{DEM}. 
\begin{prop}\label{multiplicity}
Let $(R,\m)$ be a $d$-dimensional F-finite regular local ring of characteristic $p>0$ and $\a \subset R$ an $\m$-primary ideal. 
Denote by $e(\a)$ the multiplicity of $\a$. Then
$$e(\a) \ge \left(\frac{d}{\mathrm{c}(\a)}\right)^d.$$
Moreover, the equality holds if and only if there exists an integer $n>0$ such that the integral closure $\overline{\a}$ of $\a$ is equal to $\m^n$. In this case, we have $n=d/\mathrm{c}(\a)$. 
\end{prop}
\begin{proof}
The proof is essentially the same as in \cite[Theorem 1.2, 1.4]{DEM} (cf. \cite[Theorem 3.9]{Mu2}). 
We may reduce to the case where $R$ is the localization of a $d$-dimensional polynomial ring $k[X_1, \dots, X_d]$ at the origin. 

First suppose that $\a \subset R$ is generated by monomials of $X_1, \dots, X_d$.
By \cite[Corollary 3.5]{T}, 
$$\mathrm{c}(\a)=\sup\{t \in \R_{\ge 0} \mid \tau(R,\a^t)=R \}, $$
where $\tau(R,\a^t)$ denotes an ideal associated to $\a$ with exponent $t \ge 0$ which is a generalization of test ideals introduced by Hara and Yoshida (see \cite{HY} for its definition and basic properties). 
By virtue of \cite[Theorem 4.8]{HY}, the ideal $\tau(R,\a^t)$ is characterized by the Newton polytope $t \cdot P(\a)$ associated to $\a$ and $t$. 
Hence $$\mathrm{c}(\a)=\sup\{t>0 \mid \mbox{\boldmath $e$} \in t \cdot P(\a)\},$$
where $\mbox{\boldmath $e$}=(1, \dots, 1) \in \R_{\ge 0}^d$. 
Thus, by the same argument as in \cite[Theorem 1.1]{DEM}, the assertion follows in the monomial case. 

Now we consider an arbitrary $\m$-primary ideal $\a \subset R$. 
Fixing a monomial order $\lambda$ on $k[X_1, \dots, X_d]$, we have a flat family $\{\a_s\}_{s \in k}$ such that $R/\a_s \cong R/\a$ for all $s \ne 0$ and $\a_0=\mathrm{in}_{\lambda}(\a)$. 
Since this is a flat deformation, $\mathrm{in}_{\lambda}(\a)$ is also $\m$-primary and $l(R/\a)=l(R/\mathrm{in}_{\lambda}(\a))$. 

\begin{cl}
$$\mathrm{c}(\a) \ge \mathrm{c}(\mathrm{in}_{\lambda}(\a)).$$
\end{cl}
\begin{proof}[Proof of Claim]
It suffices to show that if the pair $(R,\mathrm{in}_{\lambda}(\a)^t)$ is F-pure for some $t>0$, then the pair $(R,\a^t)$ is also F-pure. 
Suppose that $g_1, \dots, g_m$ is a Gr\"obner basis for $\a$ with respect to $>_\lambda$.
If the pair $(R,\mathrm{in}_{\lambda}(\a)^t)$ is F-pure, then by Remark \ref{split} and Lemma \ref{Fedder}, for all large $q=p^e \gg 0$, there exist nonnegative integers $l_1, \dots, l_m$ with $l_1+ \dots+l_m=\lfloor t(q-1) \rfloor$ such that $\mathrm{in}_{\lambda}(g_1)^{l_1} \cdots \mathrm{in}_{\lambda}(g_m)^{l_m} \notin \m^{[q]}$. 
Then clearly $g_1^{l_1} \cdots g_m^{l_m} \notin \m^{[q]}$, which implies the F-purity of the pair $(R, \a^t)$ by Lemma \ref{Fedder} again. 
\end{proof}
By the above claim, we have the required inequality from the case where $\a$ is a monomial ideal. 
As for the equality, our argument is the same as in \cite[Theorem 1.4]{DEM}. 
\end{proof}

\begin{rmk}
Thanks to Proposition \ref{klt}, the propositions in this section give proofs of the corresponding statements for lc threshold, which are much shorter than the original ones. 
We believe that by way of F-pure threshold, we can give algebraic proofs of  many other properties concerning lc threshold. 

\end{rmk}

\end{document}